\documentclass[10pt]{article}

\usepackage[a4paper, margin=2.25cm, top=3cm ]{geometry}

\usepackage[all]{xy}
\usepackage[nottoc]{tocbibind}
\usepackage{amsmath,amsthm,amsfonts,longtable,verbatim,tikz-cd,multicol,amssymb,wasysym,setspace,graphicx,titlesec,imakeidx,dynkin-diagrams,etoolbox}
\usepackage[inline]{enumitem}
\usepackage[colorlinks,unicode,linktocpage=true]{hyperref}
\titleformat{\subsubsection}[runin]{\normalfont\bfseries}{\thesubsubsection.}{5pt}{}{}
\titlespacing{\subsubsection}{0pt}{5pt}{10pt}

\newcommand{\bb}{\medbreak}
\newcommand{\nt}{\noindent}

\newcommand{\N}{\mathbb{N}}
\newcommand{\Cc }{\mathbb{C}}
\newcommand{\rt}{\xrightarrow{}}

\newcommand{\id}{\text{id}}

\newcommand{\End}{\text{End}}

\newcommand{\im}{\text{im}}
\newcommand{\spn}{\text{span}}

\newcommand{\Mod}{\text{-Mod}}

\newcommand{\cR}{\mathcal{R}}

\def\row#1/#2!{#1_{\IfStrEq{#2}{}{n}{#2}} & \dynkin{#1}{#2}\\}

\newcommand\blfootnote[1]{
  \begingroup
  \renewcommand\thefootnote{}\footnote{#1}
  \addtocounter{footnote}{-1}
  \endgroup
}

\newcommand{\define}[1]{\textbf{#1}\index{#1}}

\DeclareMathSymbol{\nmid}{\mathrel}{AMSb}{"2D}

\renewcommand{\thefootnote}{\fnsymbol{footnote}}

\newtheorem{lemma}{Lemma}[section]
\newtheorem{theorem}[lemma]{Theorem}

\theoremstyle{definition}
\newtheorem{definition}[lemma]{Definition}
\newtheorem{example}[lemma]{Example}

\newtheorem{remark}[lemma]{Remark}

\makeindex[intoc]

\titleformat{\subsection}[runin]{\normalfont\bfseries}{\thesubsection.}{1em}{}{}
\titlespacing{\subsection}{0pt}{5pt}{10pt}

\title{Drinfeld twists of Koszul algebras}
\author{Edward Jones-Healey}
\date{}

\counterwithin*{section}{part}

\begin{document}

\maketitle

\begin{abstract} Given a Hopf algebra $H$ and a counital $2$-cocycle $\mu$ on $H$, Drinfeld introduced a notion of twist which deforms an $H$-module algebra $A$ into a new algebra $A_\mu$. We show that when $A$ is a quadratic algebra, and $H$ acts on $A$ by degree-preserving endomorphisms, then the twist $A_\mu$ is also quadratic. Furthermore, if $A$ is a Koszul algebra, then $A_\mu$ is a Koszul algebra. As an application, we prove that the twist of the $q$-quantum plane by the quasitriangular structure of the quantum enveloping algebra $U_q(\mathfrak{sl}_2)$ is a quadratic algebra equal to the $q^{-1}$-quantum plane.
\end{abstract}

\tableofcontents

\section{Introduction}
\nt Let $k$ be an arbitrary field, and consider all objects to be linear over $k$, with $\otimes=\otimes_k$. If $H$ is a Hopf algebra over $k$, then let $\triangle:H\rt H\otimes H$ denote the coproduct and $\epsilon:H\rt k$ the counit. For an $H$-module $A$ we denote the action of $h\in H$ on $a\in A$ as $h\rhd a$.  An $H$-module algebra $A$ is an algebra for which the product map $m:A\otimes A\rt A$ is an $H$-module homomorphism, so $h\rhd m(a\otimes b)=m(\triangle(h)\rhd a\otimes b)$, and the unit $1_A$ satisfies $h\rhd 1_A=\epsilon(h)1_A$.\blfootnote{I wish to thank Yuri Bazlov for providing numerous helpful ideas and comments. Research for this paper was supported by the Engineering and Physical Sciences Research Council, UK.}\bb

\nt In \cite{drinfeld_quantum_groups}, Drinfeld introduced a notion of twist for Hopf algebras and their module algebras. We will briefly recall here how to perform this twist using the terminology of Majid \cite[Section 2.3]{alma998944944401631}. See also Etingof and Gelaki \cite[Section 5.14]{etingof2016tensor} for more. 
\begin{definition}\cite[Example 2.3.1]{alma998944944401631} A \define{$2$-cocycle} of the Hopf algebra $H$ is an invertible $\mu\in H\otimes H$ satisfying
\begin{equation}\label{2_cocycle_def}
(\mu\otimes 1)\cdot (\triangle\otimes \id)(\mu) =(1\otimes \mu)\cdot (\id\otimes \triangle)(\mu)
\end{equation}
A $2$-cocycle is said to be \define{counital} if
\begin{equation}\label{counital_def}
(\epsilon\otimes \id)(\mu)=1 =(\id\otimes \epsilon)
\end{equation}
Counital $2$-cocycles are also referred to as \textit{bialgebra twists} in Etingof and Gelaki \cite[Definition 5.14.1]{etingof2016tensor}, and \textit{twisting elements} in Giaquinto and Zhang \cite[Definition 1.2]{giaquinto1998bialgebra}.
\end{definition}
\nt The Drinfeld twist of $H$ by $\mu$ is defined to be the Hopf algebra $H_\mu$ which shares the same algebra structure and counit as $H$, and has coproduct $\triangle_\mu:H_\mu\rt H_\mu\otimes H_\mu,\ h\mapsto \mu\cdot \triangle(h)\cdot \mu^{-1}$. We shall refer to $H_\mu$ as a ``Hopf algebra twist" in order to differentiate it from the following, which is also called a Drinfeld twist.
\begin{definition} The \define{Drinfeld twist} of an $H$-module algebra $(A,m)$ by $\mu$ is the $H_\mu$-module algebra $(A_\mu,m_\mu)$ where $A_\mu=A$ as $k$-vector spaces, with the product map $m_\mu:A_\mu\otimes A_\mu\rt A_\mu$ given by $a\otimes b\mapsto m(\mu^{-1}\rhd a\otimes b)\ \forall a, b\in A_\mu$.
\end{definition} 

\nt In Section \ref{quadratic_sec} we consider the case when $A$ is a quadratic $H$-module algebra, and show that when the action of $H$ on $A$ is degree-preserving, then any Drinfeld twist of $A$ by a counital $2$-cocycle of $H$ is also a quadratic algebra. Using this, we prove in Section \ref{koszulity_sec} that if $A$ is also a Koszul algebra, then the Drinfeld twist $A_\mu$ must be Koszul too. Finally in Section \ref{examples_sec} we give several examples which apply these results.

\section{Quadratic algebras}\label{quadratic_sec}

Suppose $A$ is a connected $\N$-graded $k$-algebra, meaning $A=\bigoplus_{i\in \N}A_i$ with $A_i\cdot A_j\subseteq A_{i+j}$ and $A_0=k$. We will assume throughout this paper that $A$ is locally finite-dimensional, meaning that each grading component $A_i$ is finite-dimensional. Now $A$ is also a \define{quadratic algebra} if it is generated by its degree $1$ elements $V:=A_1$, and has quadratic relations $R\subseteq V\otimes V$, so that $A=T(V)/(R)$, where $T(V)$ is the tensor algebra over $V$ and $(R)$ is the $2$-sided ideal generated by $R$ inside $T(V)$.\bb

\nt Suppose $A=T(V)/(R)$ is both a quadratic algebra, and an $H$-module algebra, for some Hopf algebra $H$ acting by degree-preserving endomorphisms (i.e. $h\rhd a\in A_i\ \forall h\in H,a\in A_i$). Then, it easy to see that the degree $1$ subspace $V:=A_1$ is an $H$-submodule of $A$. Additionally, this action of $H$ on $V$ extends naturally by means of the coproduct of $H$ to make $T(V)$ an $H$-module algebra.\bb 

\nt Now a fact that will be very useful throughout the paper is that the space of relations $R$ is an $H$-submodule of $T(V)$. This can be seen as follows. Note that a quotient of an $H$-module algebra by an ideal can only be an $H$-module algebra too if the ideal is also an $H$-submodule. Since the $H$-module algebra structure of $A$ arises from quotienting the $H$-module algebra $T(V)$ by the ideal $(R)$, we see that $(R)$ must also be an $H$-submodule of $T(V)$. Additionally, since $R$ is precisely the subspace of degree $2$ elements in $(R)$, it must be fixed under the degree-preserving endomorphisms of $H$. And so $R$ is an $H$-submodule of $T(V)$.\bb

\nt We now give our first main result.

\begin{theorem}\label{quadratic_h_mod_alg_thm} Let $A=T(V)/(R)$ be a quadratic $H$-module algebra, where $H$ is a Hopf algebra acting by degree-preserving endomorphisms. If $\mu$ is a counital $2$-cocycle of $H$, then the Drinfeld twist $A_\mu$ is a quadratic algebra of the form $T(V)/(R_\mu)$ where $R_\mu:=\{\mu\rhd r\ |\ r\in R\}$.
\begin{proof}
Firstly let us show $A_\mu$ is a connected $\N$-graded algebra under the same grading as on $A$. Since $A_\mu=A$ as vector spaces, $A_\mu$ forms an $\N$-graded vector space under the grading of $A$. Let $m$ and $m_\mu$ denote the products on $A$ and $A_\mu$ respectively. Now $A_\mu$ is a graded algebra if $m_\mu(v\otimes w)\in A_{i+j}\ \forall v\in A_i,w\in A_j$. But this follows since $m_\mu:=m(\mu^{-1}\rhd\ )$, and the action of $H$ is degree-preserving so $\mu^{-1}\rhd v\otimes w\in A_i\otimes A_j\ \forall v\in A_i,w\in A_j$. $A_\mu$ is also connected as $(A_\mu)_0=A_0=k$.\bb

\nt Note, since $A=A_\mu$ as $\N$-graded vector spaces, $V:=A_1$ is additionally the subspace of degree $1$ elements of $A_\mu$. Also, as the action of $H$ is degree-preserving, we see that $R_\mu$ is a subspace of $V\otimes V$. Therefore, by showing $A_\mu=T(V)/(R_\mu)$, we will have proven that $A_\mu$ is a quadratic algebra, as required.\bb

\nt First we prove $A_\mu$ is generated by $V$ by induction on the degree of elements of $A_\mu$. Consider a degree $2$ element $x\in A_2$. Since $A$ is generated by $V$, $x=\sum_i m(v_i\otimes w_i)$ for some $v_i,w_i\in V$. Therefore,
$$x=\sum_i m(v_i\otimes w_i)=\sum_i m(\mu^{-1}\rhd (\mu\rhd v_i\otimes w_i))=\sum_i m_\mu(\mu\rhd v_i\otimes w_i)$$
Since the action of $H$ is degree preserving, $\mu\rhd v_i\otimes w_i\in V\otimes V$, and so $x$ can be expressed as a linear combination of $m_\mu$-products of elements of $V$, as required.\bb 

\nt Now let $k\geq2$, and suppose every element of degree $k$ in $A_\mu$ can be expressed as a linear combination of $m_\mu$-products of elements of $V$. We show that this implies the elements of degree $k+1$ have a similar such decomposition. Consider $x\in A_{k+1}$. Since $A$ is generated by $V$, we may decompose $x$ as a linear combination of terms which are the $m$-products of $k+1$-elements in $V$, i.e. is of the form $m(v_1\otimes m(v_2\otimes m(\dots)))$, for some $v_1,\dots,v_{k+1}\in V$. Suppose $\mu=\sum \mu_1\otimes \mu_2$ for some $\mu_1,\mu_2\in H$. Then  
\begin{equation}\label{a_mu_generated_by_v}
m(v_1\otimes m(v_2\otimes m(\dots )))=\sum m_\mu([\mu_1\rhd v_1]\otimes [\mu_2\rhd m(v_2\otimes m(\dots ))])
\end{equation}
Now $\mu_2\rhd m(v_2\otimes m(\dots ))\in A_k$, and so by the induction hypothesis this can be decomposed into a linear combination of $m_\mu$-products of $k$-elements of $V$. Therefore $m(v_1\otimes m(v_2\otimes m(\dots )))$ can be decomposed using $m_\mu$-products of $k+1$-elements, and this implies the same holds for $x$ too. This concludes the induction argument proving $A_\mu$ is generated by $V$.\bb 

\nt There now exists a natural surjective algebra homomorphism $\phi:T(V)\rt A_\mu$. We show $\ker(\phi)=(R_\mu)$ to complete the proof that $A_\mu$ is a quadratic algebra. Note that for the natural embedding $i:V\hookrightarrow A$, we have $R=\ker(m\circ (i\otimes i))$. As $i$ is an $H$-module homomorphism,  $i\otimes i$ is an $H\otimes H$-module homomorphism, and we see $m_\mu \circ (i\otimes i)=m\circ (i\otimes i)(\mu^{-1}\rhd\ -\ )$ . So $r\in \ker(m_\mu\circ (i\otimes i))$ iff $\mu^{-1}\rhd r\in R$, which holds iff $r\in R_\mu$. So $(R_\mu)\subseteq \ker(\phi)$, and therefore $A_\mu=T(V)/\ker(\phi)$ is seen to be a quotient of $T(V)/(R_\mu)$.\bb 

\nt To show that $A_\mu$ is equal to $T(V)/(R_\mu)$, we will apply a dimensional argument. Note that all dimensions in the following are finite due to our assumption that $A$ is locally finite-dimensional. Now, since $A_\mu$ is a quotient of $T(V)/(R_\mu)$, the dimension of each grading component of $A_\mu$ must be less than or equal to the dimension of the same degree component of $T(V)/(R_\mu)$, i.e.
\begin{equation}\label{dimensions_in_each_degree_1}
\dim(A_\mu)_i\leq \dim(T(V)/(R_\mu))_i
\end{equation}
Next we will consider the Drinfeld twist of $T(V)/(R_\mu)$ by the cocycle $\mu^{-1}$, which leads us to another dimensional inequality (see \eqref{dimensions_in_each_degree_2}). First we must check that we can indeed take such a twist. Note that this twist is over the Hopf algebra $H_\mu$, rather than the Hopf algebra $H$ that we have used so far. It is a simple exercise to check $\mu^{-1}$ is a counital $2$-cocycle of the Hopf algebra twist $H_\mu$, and we show next that $T(V)/(R_\mu)$ is an $H_\mu$-module algebra.\bb 

\nt Since $H_\mu=H$ as algebras, the $H$-action on $V$ defines an $H_\mu$-action on $V$. This extends, by means of the twisted coproduct $\triangle_\mu$, to make $T(V)$ an $H_\mu$-module algebra. Now $R_\mu$ is an $H_\mu$-submodule of $T(V)$ since, if $h\in H_\mu$ and $r\in R_\mu$, where $r=\mu\rhd r'$ for some $r'\in R$, we have
$$h\rhd r=(\mu\cdot \triangle(h)\cdot \mu^{-1})\rhd (\mu\rhd r')\in R_\mu$$
where we use the fact that $\triangle(h)\rhd r'\in R$, since we proved above that $R$ is an $H$-submodule of the $H$-module algebra $T(V)$. So we have established $T(V)$ is an $H_\mu$-module algebra, and $R_\mu$ is an $H_\mu$-submodule, and therefore $T(V)/(R_\mu)$ is an $H_\mu$-module algebra as we required.\bb

\nt We can now consider the Drinfeld twist of $T(V)/(R_\mu)$ by $\mu^{-1}$. In direct analogy to the argument used for $A_\mu$, one may show that $V$ generates $(T(V)/(R_\mu))_{\mu^{-1}}$, i.e. use the fact $T(V)/(R_\mu)$ is generated by $V$ to express elements of $(T(V)/(R_\mu))_{\mu^{-1}}$ as linear combinations of $T(V)/(R_\mu)$-products of elements of $V$. Then one rewrites each $T(V)/(R_\mu)$-product as a linear combination of $(T(V)/(R_\mu))_{\mu^{-1}}$-products.\bb 

\nt We therefore have a surjective map $\phi':T(V)\rt (T(V)/(R_\mu))_{\mu^{-1}}$. It is easy to show that $(R)\subseteq \ker(\phi')$, and so we find $(T(V)/(R_\mu))_{\mu^{-1}}$ is a quotient of $A$. This implies the following inequality in the dimensions of the grading components of degree $i$,
\begin{equation}\label{dimensions_in_each_degree_2}
\dim((T(V)/(R_\mu)_{\mu^{-1}})_i\leq \dim(A)_i
\end{equation}
But as discussed at the start of the proof, twisting preserves the grading on algebras, so
\begin{equation}\label{dimensions_in_each_degree_3}
\dim(A)_i=\dim(A_\mu)_i,\hspace{2em}\dim(T(V)/(R_\mu))_i=\dim((T(V)/(R_\mu)_{\mu^{-1}})_i
\end{equation}
Combining (in)equalities \eqref{dimensions_in_each_degree_1}, \eqref{dimensions_in_each_degree_2} and \eqref{dimensions_in_each_degree_3} we find that $\dim(A_\mu)_i=\dim(T(V)/(R_\mu))_i$. But since $A_\mu$ was a quotient of $T(V)/(R_\mu)$, we find that the algebras must be equal.
\end{proof}
\end{theorem}

\section{Koszul algebras}\label{koszulity_sec}
\subsection{Statement of the main theorem.}
Let $A$ be a connected, $\N$-graded, and locally finite-dimensional, $k$-algebra. Define $A_+:=\bigoplus_{i>0}A_i$, so $A/A_+\cong k$, and we call the corresponding quotient map $\epsilon':A\rt k$ the augmentation map. Now $k$ is an $A$-module via $\epsilon'$, and $A$ is called a \define{Koszul algebra} if there is a linear graded free resolution of $k$ as an $A$-module (see Witherspoon \cite[Definition 3.4.3]{alma992981689925201631}). If $A$ is Koszul, then it is a standard fact that it is also quadratic, so $A=T(V)/(R)$ for $V=A_1$ and $R$ a subspace of $V\otimes V$.\bb

\nt The next result establishes that Koszulity is preserved under Drinfeld twists. This generalises to arbitrary Hopf algebras a result of Davies \cite[Proposition 4.25]{doi:10.1080/00927872.2016.1178271} who proved this for the case when the Hopf algebra $H$ is the group algebra of a finite abelian group. The rest of the section is dedicated to proving this result.

\begin{theorem}\label{koszul_main_thm} Let $A=T(V)/(R)$ be a Koszul $H$-module algebra, where $H$ is a Hopf algebra acting by degree-preserving endomorphisms. If $\mu$ is a counital $2$-cocycle of $H$, then the Drinfeld twist $A_\mu$ is a Koszul algebra given by $T(V)/(R_\mu)$ where $R_\mu:=\{\mu\rhd r\ |\ r\in R\}$. 
\end{theorem}

\subsection{Plan for the proof.}
It follows immediately from Theorem \ref{quadratic_h_mod_alg_thm} that $A_\mu$ is a quadratic algebra of the form $T(V)/(R_\mu)$. Using this we can construct a complex $K_\bullet(A_\mu)$, which, by Witherspoon \cite[Theorem 3.4.6]{alma992981689925201631}, is a resolution of $k$ as an $A_\mu$-module precisely when $A_\mu$ is a Koszul algebra. We therefore work to show $K_\bullet(A_\mu)$ is indeed a resolution to complete the proof. We start by considering the Koszul resolution $K_\bullet(A)$ of $A$, and twist this resolution using a functor of Giaquinto and Zhang. This produces a new resolution of $k$ as an $A_\mu$-module, which we denote $K_\bullet(A)_\mu$. We then construct an isomorphism of complexes between $K_\bullet(A)_\mu$ and $K_\bullet(A_\mu)$, the existence of which implies $K_\bullet(A_\mu)$ is a resolution $k$ as an $A_\mu$-module, as required.

\subsection{The Koszul resolution of \texorpdfstring{$A$}{A}.}
Let us start by defining the Koszul resolution of $k$ as an $A$-module. It is given by $K_\bullet(A)=\bigoplus_{n\geq 0} K_n(A)$, where
\begin{equation}\label{koszul_complex}
K_0(A)=A,\hspace{2em} K_1(A)=A\otimes V,\hspace{2em} K_n(A)=A\otimes \bigcap_{i+j=n-2}(V^{\otimes i}\otimes R\otimes V^{\otimes j})\ \text{for }n\geq 2
\end{equation}
The differentials $d_n$ are induced by the canonical embedding of $K_\bullet(A)$ into the Bar resolution $B_\bullet(A)$ of $k$, where $B_n(A)=A^{\otimes (n+1)}$ and 
\begin{equation}\label{differential}
d_n(a_0\otimes \dots \otimes a_n)=(-1)^n\epsilon'(a_n)a_0\otimes \dots \otimes a_{n-1}+\sum_{i=0}^{n-1}(-1)^i a_0\otimes \dots \otimes a_i a_{i+1}\otimes \dots \otimes a_n
\end{equation}
for $a_0,\dots ,a_n\in A$ and $\epsilon'$ is the augmentation map of $A$. $K_n(A)$ is a left $A$-module under multiplication on the leftmost tensor leg. It is also a left $H$-module by restricting the natural action of $H$ on $B_n(A)=A^{\otimes (n+1)}$ (using the coproduct $\triangle$ of $H$) onto $K_n(A)$. $K_n(A)$ is closed under this $H$-action since we showed above that $V$ and $R$ are $H$-submodules of $T(V)$, so can also be seen as $H$-submodules of $A$ and $A\otimes A$ respectively. Therefore $K_n(A)$ is an intersection of $H$-submodules of $B_n(A)$, so is an $H$-submodule itself.

\subsection{The Giaquinto and Zhang twisting functor.} Let $A\Mod$ be the category of all left $A$-modules. Take the category $(H,A)\Mod$ to be the subcategory of $A\Mod$ whose objects $M$ are also left $H$-modules such that the following holds for all $m\in M, h\in H, a\in A$,
\begin{equation}\label{H_A_module}
h\rhd (a\rhd m)=(h_{(1)}\rhd a)\rhd (h_{(2)}\rhd m)
\end{equation}
where $\triangle(h)=h_{(1)}\otimes h_{(2)}$. The morphisms are those of $A\Mod$ that are also $H$-module homomorphisms.\bb

\nt So far, Drinfeld twists have provided a mechanism for deforming the Hopf algebra $H$ and the $H$-module algebra $A$. Giaquinto and Zhang \cite[Theorem 1.7]{giaquinto1998bialgebra} extends this by defining a way to twist a module in $(H,A)\Mod$ into a module in $(H_\mu,A_\mu)\Mod$. This new twist defines an equivalence of categories $(H,A)\Mod\rt (H_\mu,A_\mu)\Mod$. We show next that the Koszul resolution $K_\bullet(A)$ can be defined within $(H,A)\Mod$, and describe its image under this twisting functor of Giaquinto and Zhang.\bb

\nt Let us first check that \eqref{H_A_module} holds on $K_n(A),\ \forall n\geq 0$. For $K_0(A)=A$, let $h\in H,a,b\in A$, then 
$$h\rhd (a\rhd b)=h\rhd (a\cdot b)=(h_{(1)}\rhd a)\cdot (h_{(2)}\rhd b)=(h_{(1)}\rhd a)\rhd (h_{(2)}\rhd b)$$
as required. For $n\geq 1$, let $m\in K_n(A)\subseteq A^{\otimes (n+1)}$. Suppose $m=a'\otimes m'$ for some $a'\in A,\ m'\in A^{\otimes n}$. Then,
\begin{align*}
h\rhd (a\rhd m)=h\rhd (a\cdot a'\otimes m') & =(h_{(1)}\rhd a\otimes a')\otimes (h_{(2)}\rhd m')\\
& = (h_{(1)(1)}\rhd a)\cdot (h_{(1)(2)}\rhd a') \otimes (h_{(2)}\rhd m')\\
& = (h_{(1)}\rhd a)\cdot (h_{(2)(1)}\rhd a') \otimes (h_{(2)(2)}\rhd m')\\
& =(h_{(1)}\rhd a)\rhd (h_{(2)}\rhd m)
\end{align*} 
So \eqref{H_A_module} is satisfied and $K_n(A)$ is an object in $(H,A)\Mod$. It is standard that the differentials $d_n$ of the Koszul resolution are $A$-module homomorphisms, but we note that they are also $H$-module homomorphisms. Indeed, on inspecting \eqref{differential}, we see that as a map,
\begin{equation}\label{differential_map}
d_n=\sum^{n-1}_{i=0}(-1)^i \id^{\otimes i}\otimes m\otimes \id^{\otimes n-i-1}+(-1)^n \id^{\otimes n}\otimes \epsilon'
\end{equation}
where $m$ denotes the product map of $A$ and $\epsilon'$ is the augmentation map of $A$. Since $m$ and $\epsilon'$ are $H$-module homomorphisms, it follows that $d_n$ also is. Therefore the whole Koszul resolution $K_\bullet(A)$ can be defined within the category $(H,A)\Mod$.\bb

\nt We now apply the functor of Giaquinto and Zhang \cite[Theorem 1.7]{giaquinto1998bialgebra}. Under this functor the $A$-modules $K_n(A)$ are twisted into $A_\mu$-modules $K_n(A)_\mu$ in the following way: let $K_n(A)_\mu= K_n(A)$ as $k$-vector spaces, and equip $K_n(A)_\mu$ with the action $$\rhd_\mu:A_\mu\otimes K_n(A)_\mu\rt K_n(A)_\mu,\ \rhd_\mu=\rhd (\mu^{-1}\rhd \ -\ )$$
Here we use the fact that $A$ and $K_n(A)$ are $H$-modules, so $A\otimes K_n(A)$ is naturally an $H\otimes H$-module. Therefore $\mu^{-1}\rhd\ $ defines an endomorphism of $A\otimes K_n(A)$, and also of $A_\mu\otimes K_n(A)_\mu$, since these are equal as $k$-vector spaces.\bb

\nt We must also apply the functor to the differentials in the Koszul resolution, but Giaquinto and Zhang do not explicitly describe how their functor behaves on morphisms. For completeness we prove here that we can take it as mapping each morphism in $(H,A)\Mod$ to itself. Suppose $V,W\in (H,A)\Mod$ with the actions of $A$ on $V$ and $W$ denoted by $\rhd^V$ and $\rhd^W$ respectively. Let $\phi:V\rt W$ be an $A$-module homomorphism, so $\rhd^W\circ (\id_A\otimes \phi)=\phi\circ \rhd^V$. Then $\phi$ is also an $A_\mu$-module homomorphism $V_\mu\rt W_\mu$. Indeed 
$$\rhd^W_\mu\circ (\id_A\otimes \phi)=\rhd^W\circ(\id_A\otimes \phi) (\mu^{-1}\rhd \ -\ )=(\phi\circ \rhd^V)(\mu^{-1}\rhd \ -\ )=\phi\circ \rhd^V_\mu$$
where we use the fact $\id_A\otimes \phi$ is an $H$-module homomorphism in the first equality. $\phi$ is also an $H_\mu$-module homomorphism $V_\mu\rt W_\mu$, since $H_\mu$ acts on $V_\mu$ and $W_\mu$ in exactly the same way that $H$ acts on $V$ and $W$ respectively. Therefore $\phi$ can also be viewed as a morphism in $(H_\mu,A_\mu)\Mod$ from $V_\mu$ to $W_\mu$, and it makes sense to let the functor of Giaquinto and Zhang act identically on morphisms. Therefore each differential $d_n$ of the Koszul resolution $K_\bullet(A)$ will be sent by the functor of Giaquinto and Zhang to itself.\bb

\nt To summarise, the result of applying the twisting functor to $K_\bullet(A)$ is a complex $K_\bullet(A)_\mu$ of $A_\mu$-modules which share the same underlying vector spaces, and differentials, as those on $K_\bullet(A)$. Since $K_\bullet(A)$ is a resolution of $k$ as an $A$-module, we see that the twisted complex $K_\bullet(A)_\mu$ must be a resolution of $k$ as an $A_\mu$-module.

\subsection{The Koszul complex of \texorpdfstring{$A_\mu$}{A}.} Using Theorem \ref{quadratic_h_mod_alg_thm} we know that $A_\mu$ is a quadratic algebra of the form $T(V)/(R_\mu)$. We can therefore construct a complex $K_\bullet(A_\mu)$ completely analogously to the Koszul resolution for $A$ given in \eqref{koszul_complex}, and we describe this explicitly next. Let $K_\bullet(A_\mu)=\bigoplus_{n\geq 0} K_n(A_\mu)$ where
\begin{equation*}
   K_0(A_\mu)=A_\mu,\hspace{2em} K_1(A_\mu)=A_\mu\otimes V,\hspace{2em} K_n(A_\mu)=A_\mu\otimes \bigcap_{i+j=n-2}(V^{\otimes i}\otimes R_\mu\otimes V^{\otimes j})\ \text{for }n\geq 2
\end{equation*}
This is a complex of $A_\mu$-modules, where the differentials are inherited from the Bar resolution $B_\bullet(A_\mu)$, where $B_n(A_\mu)=(A_\mu)^{\otimes (n+1)}$ and
\begin{equation}\label{differential_mu}
d^\mu_n=\sum^{n-1}_{i=0}(-1)^i \id^{\otimes i}\otimes m_\mu\otimes \id^{\otimes n-i-1}+(-1)^n \id^{\otimes n}\otimes \epsilon'
\end{equation}
where we use the same augmentation map $\epsilon'$ as on $A$, since $A_\mu$ and $A$ have the same grading.\bb

\nt Note that by Witherspoon \cite[Theorem 3.4.6]{alma992981689925201631}, $K_\bullet(A_\mu)$ is a resolution of $k$ as an $A_\mu$-module precisely when $A_\mu$ is a Koszul algebra. We use this to prove the theorem, showing that $K_\bullet(A_\mu)$ is indeed a resolution. To do so we construct an isomorphism of complexes from $K_\bullet(A)_\mu$ to $K_\bullet(A_\mu)$, and use the fact established in the last section that $K_\bullet(A)_\mu$ is a resolution of $k$ as an $A_\mu$-module. To construct this isomorphism of complexes we first define a sequence of elements which turn out to generalise counital $2$-cocycles.

\subsection{Higher counital \texorpdfstring{$2$}{2}-cocycles.} First we introduce some notation: for $i\in \{1,\dots,n\}$, let $\triangle^{(n)}_i:H^{\otimes n}\rt H^{\otimes n+1},\ \triangle^{(n)}_i=\id\otimes \dots \otimes \triangle\otimes \dots \otimes \id$, where $\triangle$ appears in the $i$-th position. By convention also let $\triangle^{(n)}_0(-)=1\otimes (-)$ and $\triangle^{(n)}_{n+1}(-)=(-)\otimes 1$.\bb

\nt Recall that our counital $2$-cocycle $\mu$ satisfies the $2$-cocycle equation in \eqref{2_cocycle_def}. This equation can be rewritten in the above notation as: 
\begin{equation}\label{2_cocycle_eqn}
\triangle^{(2)}_3(\mu)\cdot \triangle^{(2)}_1(\mu)=\triangle^{(2)}_0(\mu)\cdot \triangle^{(2)}_2(\mu)
\end{equation}
It turns out that $\mu$ can be viewed as just one step in a sequence of ``higher counital $2$-cocycles", where the next two lemmas justify this terminology. For $n \geq 0$, let $f_n\in H^{\otimes (n+1)}$ be defined as $f_0 = 1$, $f_1 = \mu$ and for $n\geq 2$,
\begin{equation}\label{f_n}
f_n = (\mu\otimes 1^{\otimes n-1})\cdot \prod_{i=0}^{n-2}\big(\triangle^{(n)}_1 \circ \triangle^{(n-1)}_1 \circ \dots \circ \triangle^{(n-i)}_1\big) (\mu\otimes 1^{\otimes n-i-2})
\end{equation}
where $\circ$ denotes composition, and all products are taken in the algebra $H^{\otimes (n+1)}$. Notice that $f_2=(\mu\otimes 1)\cdot (\triangle\otimes \id)(\mu)$, which is just the left hand side of \eqref{2_cocycle_eqn}. The following lemma can be viewed as saying that each of the elements $f_n$ satisfies a kind of ``higher" version of the $2$-cocycle equation \eqref{2_cocycle_eqn}:
\begin{lemma}\label{higher_2_cocycle} For all $n\in \N$ and $i\in \{1,\dots,n\}$, $f_n=(1^{\otimes i-1}\otimes \mu\otimes 1^{\otimes n-i})\cdot\triangle^{(n)}_i(f_{n-1})$.
\begin{proof}
Let $f^i_n:=(1^{\otimes i-1}\otimes \mu\otimes 1^{\otimes n-i})\cdot\triangle^{(n)}_i(f_{n-1})$. Firstly it is clear that $f_1=f^1_1$. Now let us check $f_n=f^1_n$ for $n\geq 2$:
\begin{align*}
f_n^1 & = (\mu\otimes 1^{\otimes n-1})\cdot \triangle^{(n)}_1(f_{n-1}) \\
 & =  (\mu\otimes 1^{\otimes n-1})\cdot \triangle^{(n)}_1 (\mu\otimes 1^{\otimes n-2})  \cdot   \triangle^{(n)}_1 \bigg(\prod_{i=0}^{n-3}(\triangle^{(n-1)}_1 \circ \dots \circ \triangle^{(n-i-1)}_1) (\mu\otimes 1^{\otimes n-i-3})\bigg)\\
 & = (\mu\otimes 1^{\otimes n-1})\cdot \prod_{i=0}^{n-2}\big(\triangle^{(n)}_1\circ \triangle^{(n-1)}_1 \circ \dots \circ \triangle^{(i+2)}_1\big) (\mu\otimes 1^{\otimes i}) = f_n
\end{align*}
To finish proving the lemma we show that, for each $n\geq 2$, $f^i_n=f^{i+1}_n\ \forall i\in\{1,\dots,n-1\}$. We do so by performing induction on $n$. The base case $n=2$ holds since $f^1_2=f^2_2$ is precisely equation \eqref{2_cocycle_eqn}. Now suppose the hypothesis holds for $n=k-1$, and let us show it holds for $n=k$. Take some $i\in\{1,\dots,n-1\}$. We know that $f_{k-1}=f^1_{k-1}$, and now by the induction hypothesis we have 
\begin{equation}\label{f_k_1}
f_{k-1}=f^i_{k-1}=(1^{\otimes i-1}\otimes \mu\otimes 1^{\otimes k-i-1})\cdot \triangle^{(k-1)}_i(f_{k-2})
\end{equation}
Therefore,
\begin{align*}
f^{i+1}_k & =(1^{\otimes i}\otimes \mu\otimes 1^{\otimes k-i-1}) \cdot   \triangle^{(k)}_{i+1}\big( (1^{\otimes i-1}\otimes \mu\otimes 1^{\otimes k-i-1})\cdot\triangle^{(k-1)}_i(f_{k-2})\big)\\
& =(1^{\otimes i}\otimes \mu\otimes 1^{\otimes k-i-1}) \cdot  \triangle^{(k)}_{i+1}\big(1^{\otimes i-1}\otimes \mu\otimes 1^{\otimes k-i-1}\big) \cdot \triangle^{(k)}_{i+1}\big(\triangle^{(k-1)}_i(f_{k-2})\big)\\
& =\big(1^{\otimes i-1}\otimes [\triangle^{(2)}_{0}(\mu)\cdot \triangle^{(2)}_{2}(\mu) ]\otimes 1^{\otimes k-i-1}\big)\cdot \triangle^{(k)}_{i+1}\big(\triangle^{(k-1)}_i(f_{k-2})\big)\\
& = \big(1^{\otimes i-1}\otimes [\triangle^{(2)}_{3}(\mu)\cdot\triangle^{(2)}_{1}(\mu)]\otimes 1^{\otimes k-i-1}\big)\cdot\triangle^{(k)}_{i+1}\big(\triangle^{(k-1)}_i(f_{k-2})\big)\\
& =(1^{\otimes i-1}\otimes \mu\otimes 1^{\otimes k-i})\cdot \triangle^{(k)}_{i}\big(1^{\otimes i-1}\otimes \mu\otimes 1^{\otimes k-i-1}\big)\cdot \triangle^{(k)}_{i+1}\big(\triangle^{(k-1)}_i(f_{k-2})\big)
\end{align*}
where in the first equality we use the definition of $f^{i+1}_k$, and insert the expression \eqref{f_k_1} for $f_{k-1}$. Next we use the fact that $\triangle^{(k)}_{i+1}$ is an algebra homomorphism. In the 3rd equality we take the product of the first two terms in the previous line, and in the $4$-th equality we apply the $2$-cocycle equation \eqref{2_cocycle_eqn}. Finally notice that $\triangle^{(k)}_{i+1}\circ \triangle^{(k-1)}_i=\triangle^{(k)}_{i}\circ \triangle^{(k-1)}_i$ by coassociativity of $H$, and therefore the final expression above reduces to $f^i_k$, as required.
\end{proof}
\end{lemma}

\nt Recall that $\mu$ also satisfies the counital equation \eqref{counital_def}, i.e. $(\epsilon\otimes \id)(\mu)=1 =(\id\otimes \epsilon)(\mu)$, where $\epsilon$ is the counit of $H$. The following lemma proves that the elements $f_n$ satisfy a generalised notion of counitality,
\begin{lemma}\label{higher_counital_axiom}
For all $n\geq 1$ and $i\in\{0,\dots,n\}$, $(\id^{\otimes i}\otimes \epsilon\otimes \id^{\otimes n-i})(f_n)=f_{n-1}$.
\begin{proof} The $n=1$ case of this is precisely the counital axiom \eqref{counital_def}. For arbitrary $n\geq 1$ and $i\in \{0,\dots,n\}$, we can apply Lemma \ref{higher_2_cocycle} to express $f_n$ as $(1^{\otimes i-1}\otimes \mu\otimes 1^{\otimes n-i})\cdot\triangle^{(n)}_i(f_{n-1})$. Then
\begin{align*}
(\id^{\otimes i}\otimes \epsilon\otimes \id^{\otimes n-i}) & (f_n) = (\id^{\otimes i}\otimes \epsilon\otimes \id^{\otimes n-i})(1^{\otimes i-1}\otimes \mu\otimes 1^{\otimes n-i})\cdot (\id^{\otimes i}\otimes \epsilon\otimes \id^{\otimes n-i})(\triangle^{(n)}_i(f_{n-1}))
\end{align*}
Now 
$$(\id^{\otimes i}\otimes \epsilon\otimes \id^{\otimes n-i})(1^{\otimes i-1}\otimes \mu\otimes 1^{\otimes n-i})=1^{\otimes i-1}\otimes (\id\otimes \epsilon)(\mu)\otimes 1^{\otimes n-i}$$
so by the counitality of $\mu$, this reduces to $1^{\otimes n}$. Also,
$$ (\id^{\otimes i}\otimes \epsilon\otimes \id^{\otimes n-i})(\triangle^{(n)}_i(f_{n-1}))=(\id^{\otimes i-1}\otimes [(\id \otimes \epsilon)\circ \triangle]\otimes \id^{\otimes n-i})(f_{n-1})$$
and this is equal to $f_{n-1}$ since by the counit axiom for $H$, $(\id\otimes \epsilon)\circ \triangle=\id$.
\end{proof}
\end{lemma}

\subsection{Defining the isomorphism of complexes \texorpdfstring{$F_\bullet$}{PDFstring}.} 
Using the higher counital $2$-cocycles constructed in the last section we will now define an isomorphism of complexes between $K_\bullet(A)_\mu$ and $K_\bullet(A_\mu)$. From this we will deduce that $K_\bullet(A_\mu)$ is a resolution, and therefore that $A_\mu$ is a Koszul algebra.\bb

\nt Recall that the complex $K_\bullet(A)$ embeds into the Bar resolution $B_\bullet(A)$, so in particular $K_n(A)$ is a subspace of $A^{\otimes n+1}$, for all $n\geq 0$. Additionally, the twisted module $K_n(A)_\mu$ is equal to $K_n(A)$ as a $k$-vector space, and therefore it is also a subspace of $A^{\otimes n+1}$. Since $A$ is an $H$-module algebra, $A^{\otimes n+1}$ is naturally an $H^{\otimes n+1}$-module. Therefore we may consider the action of $f_n\in H^{\otimes n+1}$ on the space $K_n(A)_\mu$.\bb

\nt Define $F_\bullet:K_\bullet(A)_\mu\rt K_\bullet(A_\mu)$ to be the corresponding sequence of $k$-linear maps, i.e. for each $n\geq 0$, $F_n$ is the map given by $f_n$ acting on $K_n(A)_\mu$. We check in Section \ref{image_check_sec} that the image of $F_n$ is indeed inside $K_n(A_\mu)$. The first few maps of $F_\bullet$ are depicted in the diagram below,
\begin{equation*}\xymatrix@C=4em{
\dots \ar[r] & K_2(A_\mu) \ar[r]^{d^\mu_2}  & K_1(A_\mu) \ar[r]^{d^\mu_1}  & K_0(A_\mu) \ar[r] & 0\\
\dots \ar[r] & K_2(A)_\mu \ar[r]^{d_2} \ar[u]^{(\mu\otimes 1)\cdot (\triangle\otimes \id)(\mu)\rhd} & K_1(A)_\mu \ar[r]^{d_1} \ar[u]^{\mu\rhd} & K_0(A)_\mu \ar[r] \ar[u]^{\id} & 0
}
\end{equation*}
We will also check in Section \ref{diagram_commutes_sec} that the diagram above commutes, i.e. $d^\mu_nF_n=F_{n-1}d_n,\ \forall n\geq 1$, and in Section \ref{inverse_sec} that each $F_n$ has a $k$-linear inverse. These facts are sufficient to deduce that $K_\bullet(A_\mu)$ is a resolution of $k$ as an $A_\mu$-module, as we will explain in Section \ref{finale_sec}.

\begin{remark} 
It is also possible to show that the $k$-linear maps $F_n$ are $A_\mu$-module homomorphisms, and therefore, in combination with the above properties, $F_\bullet$ is an isomorphism of complexes. However it is not required that $F_n$ be an $A_\mu$-module homomorphism in order to prove $K_\bullet(A_\mu)$ is a resolution, so we won't include the proof of this fact here. For brevity though, we will still refer to $F_\bullet$ as being a chain map, or isomorphism of complexes, throughout the paper. 
\end{remark}

\subsection{Checking \texorpdfstring{$\im(F_n)\subseteq K_n(A_\mu)$}{Fn}.}\label{image_check_sec}
For the $n=0$ case, $F_0=\id$, and as vector spaces, $K_0(A)_\mu=A=A_\mu=K_0(A_\mu)$, so the image of $F_0$ is equal to $K_0(A_\mu)$, as required.\bb

\nt When $n=1$, we have as vector spaces $K_1(A)_\mu=A\otimes V=A_\mu\otimes V=K_1(A_\mu)$. Now $V$ was defined to be the degree $1$ component $A_1$ of $A$, and since the $H$-action on $A$ is degree-preserving, $V$ is an $H$-submodule of $A$. Therefore $A\otimes V$ is an $H\otimes H$-submodule of $A\otimes A$, so the action of $\mu$ on $K_1(A)_\mu$ will be closed. Due to the equality of vector spaces $K_1(A)_\mu=K_1(A_\mu)$, we can view the image of $\mu\rhd$ as being in $K_1(A_\mu)$.\bb

\nt For $n\geq 2$, we inspect how $f_n$ in \eqref{f_n} acts on an element $a\in K_n(A)_\mu$, and show the result lies in $K_n(A_\mu)$. Recall that, 
$$K_n(A)_\mu=A\otimes \bigcap_{i+j=n-2}(V^{\otimes i}\otimes R\otimes V^{\otimes j})\hspace{2em}K_n(A_\mu)=A_\mu\otimes \bigcap_{i+j=n-2}(V^{\otimes i}\otimes R_\mu\otimes V^{\otimes j})$$
\nt If $a\in K_n(A)_\mu$, then $a\in A\otimes V^{\otimes i} \otimes R\otimes V^{\otimes n-2-i}$ for all $i\in \{0,\dots,n-2\}$. Similarly, we can show $f_n\rhd a\in K_n(A_\mu)$ by showing $f_n\rhd a\in A\otimes V^{\otimes i} \otimes R_\mu\otimes V^{\otimes n-2-i}$ for all $i\in \{0,\dots,n-2\}$, where we have used here the fact $A_\mu=A$ as vector spaces. Fix an arbitrary $i\in \{0,\dots,n-2\}$. By Lemma \ref{higher_2_cocycle}, $f_n=(1^{\otimes i+1}\otimes \mu\otimes 1^{\otimes n-i-2}) \cdot \triangle^{(n)}_{i+2}(f_{n-1})$, and so,
$$f_n\rhd a=(1^{\otimes i+1}\otimes \mu\otimes 1^{\otimes n-i-2}) \rhd \big(\triangle^{(n)}_{i+2}(f_{n-1})\rhd a\big)$$
\nt Let us inspect where we land after $\triangle^{(n)}_{i+2}(f_{n-1})$ acts on $a$. First note $a\in A\otimes V^{\otimes i} \otimes R\otimes V^{\otimes n-2-i}$, and since $A$ and $V$ are $H$-modules, the legs of $\triangle^{(n)}_{i+2}(f_{n-1})$ hitting $A$ or $V$ will remain in these spaces. Recall also that $R$ is an $H$-submodule of $T(V)$, so $\triangle(h)\rhd r\in R,\ \forall h\in H,r\in R$. Since the $\triangle$ in $\triangle^{(n)}_{i+2}=\id\otimes \dots \otimes \triangle\otimes \dots \otimes \id$ hits the $R$ in $A\otimes V^{\otimes i} \otimes R \otimes V^{\otimes n-2-i}$, we may apply the fact that $R$ is an $H$-submodule of $T(V)$, to deduce 
$$\triangle^{(n)}_{i+2}(f_{n-1})\rhd a\in A\otimes V^{\otimes i} \otimes R\otimes V^{\otimes n-2-i}$$
Now when $1^{\otimes i+1}\otimes \mu\otimes 1^{\otimes n-i-2}$ acts on an element of $A\otimes V^{\otimes i} \otimes R\otimes V^{\otimes n-2-i}$ we clearly we land in $A\otimes V^{\otimes i} \otimes R_\mu\otimes V^{\otimes n-2-i}$.
Since $i$ was arbitrary in $\{0,\dots,n-2\}$, we find $f_n\rhd a\in K_n(A_\mu)$, and therefore the image of $F_n$ is indeed inside $K_n(A_\mu)$ as required.\bb
\subsection{\texorpdfstring{$F_\bullet$}{PDFstring} commutes with differentials.}\label{diagram_commutes_sec}
Let us now check that the maps $F_n$ make the above diagram commute, i.e. $d^\mu_nF_n=F_{n-1}d_n, \forall n\geq 1$. We inspect the result of each side of this condition on an element $a_0\otimes \dots\otimes a_n\in K_n(A)_\mu$:
\begin{align*}
F_{n-1}d_n(a_0\otimes & \dots  \otimes a_n) \\
= & f_{n-1}\rhd \bigg(\sum^{n-1}_{i=0}(-1)^i \id^{\otimes i}\otimes m\otimes \id^{\otimes n-i-1}+(-1)^n \id^{\otimes n}\otimes \epsilon'\bigg)(a_0\otimes \dots \otimes a_n)\\
 = & \sum^{n-1}_{i=0}(-1)^i (\id^{\otimes i}\otimes m\otimes \id^{\otimes n-i-1})\bigg(\triangle^{(n)}_{i+1}(f_{n-1})\rhd a_0\otimes \dots \otimes a_n \bigg)\\
& +(-1)^n f_{n-1}\rhd (\id^{\otimes n}\otimes \epsilon')\big(a_0\otimes \dots \otimes a_n\big)
\end{align*}
where in the first equality we apply the definition of the map $d_n$ given in \eqref{differential_map}. In the second equality we use the fact that $(\id^{\otimes i}\otimes m\otimes \id^{\otimes n-i-1})$ is an $H^{\otimes n}$-module homomorphism to pull the $f_{n-1}\rhd $ inside, resulting in the $\triangle^{(n)}_{i+1}(f_{n-1})\rhd $ term. Now,
\begin{align*}
d^\mu_n F_n(a_0\otimes & \dots \otimes a_n)\\
 = & \sum^{n-1}_{i=0}(-1)^i (\id^{\otimes i}\otimes m_\mu\otimes \id^{\otimes n-i-1})\bigg(f_n\rhd a_0\otimes \dots \otimes a_n \bigg) + (-1)^n (\id^{\otimes n}\otimes \epsilon')\bigg(f_n\rhd a_0\otimes \dots \otimes a_{n}\bigg)\\
 = & \sum^{n-1}_{i=0}(-1)^i (\id^{\otimes i}\otimes m\otimes \id^{\otimes n-i-1})\bigg((1^{\otimes i}\otimes \mu^{-1}\otimes 1^{\otimes n-i-1})\cdot f_n\rhd a_0\otimes \dots \otimes a_n \bigg)\\
& + (-1)^n (\id^{\otimes n}\otimes \epsilon')\bigg(f_n\rhd a_0\otimes \dots \otimes a_{n}\bigg)
\end{align*}
In the first equality we use the definition of $d_n^\mu$ in \eqref{differential_mu}. In the second we turn $m_\mu$ into $m$, which results in the $(1^{\otimes i} \otimes \mu^{-1}\otimes 1^{\otimes n-i-1})\rhd $ term.\bb

\nt By Lemma \ref{higher_2_cocycle}, $f_n=(1^{\otimes i}\otimes \mu\otimes 1^{\otimes n-i-1})\cdot\triangle^{(n)}_{i+1}(f_{n-1})$, and therefore $(1^{\otimes i}\otimes \mu^{-1}\otimes 1^{\otimes n-i-1})\cdot f_n=\triangle^{(n)}_{i+1}(f_{n-1})$. This nearly concludes the proof that $d^\mu_n F_n(a_0\otimes \dots \otimes a_n)=F_{n-1}d_n(a_0\otimes \dots  \otimes a_n)$. It only remains to show
$$(\id^{\otimes n}\otimes \epsilon')(f_n\rhd a_0\otimes \dots \otimes a_{n})=f_{n-1}\rhd\big(\id^{\otimes n}\otimes \epsilon'\big)(a_0\otimes \dots \otimes a_n)$$
\nt As $A$ is an $H$-module algebra, we have $h\rhd 1=\epsilon(h)1$ where $\epsilon$ is the counit of $H$. Therefore the degree $0$ component $A_0\cong k$ of $A$ has the structure of a trivial $H$-module, i.e. $h\rhd \lambda=\epsilon(h)\lambda,\ \forall h\in H,\lambda\in A_0$. Since the augmentation map $\epsilon'$ maps into $k$, we find $h\rhd \epsilon'(x)=\epsilon(h)\epsilon'(x)$. Using this, and the fact $\epsilon'$ is an $H$-module homomorphism, we find 
\begin{align*}
(\id^{\otimes n}\otimes \epsilon')(f_n\rhd a_0\otimes \dots \otimes a_n) & = f_n\rhd (\id^{\otimes n}\otimes \epsilon')(a_0\otimes \dots \otimes a_n)\\
& =(\id^{\otimes n}\otimes \epsilon)(f_n)\rhd (\id^{\otimes n}\otimes \epsilon')(a_0\otimes \dots \otimes a_{n})\\
& = f_{n-1}\rhd (\id^{\otimes n}\otimes \epsilon')(a_0\otimes \dots \otimes a_{n})
\end{align*}
where we apply the ``higher counitality" of $f_n$ (Lemma \ref{higher_counital_axiom}) in the final equality to say $f_{n-1}=(\id^{\otimes n}\otimes \epsilon)(f_n)$.
\subsection{The inverse chain map.}\label{inverse_sec}
Finally it is clear from the fact that $\mu$ is invertible that each map $F_n$ has a $k$-linear inverse. In particular, let $F_n^{-1}:=f_n^{-1}\rhd$, where $f_0^{-1}=1,\ f_1^{-1}=\mu^{-1}$ and for $n\geq 2$,
\begin{equation}\label{f_n_inverse}
f_n^{-1}= \bigg(\prod_{i=0}^{n-2}(\triangle^{(n)}_1\circ \triangle^{(n-1)}_1 \circ \dots \circ \triangle^{(i+2)}_1) (\mu^{-1}\otimes 1^{\otimes i})\bigg)\cdot (\mu^{-1}\otimes 1^{\otimes n-1})
\end{equation}
Although this is not required, it can be checked that $F_\bullet^{-1}$ also forms a chain map, and therefore that $F_\bullet$ is an isomorphism of complexes, as claimed.

\subsection{\texorpdfstring{$K_n(A_\mu$)}{PDFstring} is a resolution of \texorpdfstring{$k$}{PDFstring} as an \texorpdfstring{$A_\mu$}{PDFstring}-module.}\label{finale_sec}
Recall that the complex $K_\bullet(A)_\mu$, which arose from applying the functor of Giaquinto and Zhang to the Koszul resolution $K_n(A)$, is a resolution of $k$ as an $A_\mu$-module. This means that the following complex is exact,
\begin{equation}\label{augmented_complex}\xymatrix@C=4em{
\dots \ar[r] & K_2(A)_\mu \ar[r]^{d_2} & K_1(A)_\mu \ar[r]^{d_1} & K_0(A)_\mu \ar[r]^{\epsilon'} & k \ar[r] & 0
}
\end{equation}
We wish to show that $K_\bullet(A_\mu)$ is also a resolution of $k$ as an $A_\mu$-module, as this is equivalent to $A_\mu$ being a Koszul algebra. Hence we must show the following is exact,
\begin{equation*}\xymatrix@C=4em{
\dots \ar[r] & K_2(A_\mu) \ar[r]^{d^\mu_2}  & K_1(A_\mu) \ar[r]^{d^\mu_1}  & K_0(A_\mu) \ar[r]^{\epsilon'} & k \ar[r] & 0
}
\end{equation*}
In particular we wish to show $\im(d^\mu_{n+1})=\ker(d^\mu_n)$ for $n\geq 1$ and $\im(d^\mu_1)=\ker(\epsilon')$.\bb 

\nt We showed in Section \ref{diagram_commutes_sec} that $F_\bullet$ commutes with differentials, so $d^\mu_{n+1}=F_nd_{n+1}F^{-1}_{n+1},\ \forall n\geq 0$. Additionally every map $F_n$ is a $k$-linear isomorphism, and so we find $\im(d^\mu_{n+1})=F_n(\im(d_{n+1}))$. Also, for $n\geq 1$, $d^\mu_n=F_{n-1}d_nF^{-1}_n$, and we see that $\ker(d^\mu_n)=F_n(\ker(d_n))$. Therefore, for $n\geq 1$,
$$\im(d^\mu_{n+1})=F_n(\im(d_{n+1}))=F_n(\ker(d_n))=\ker(d^\mu_n)$$
where in the second equality we apply the fact \eqref{augmented_complex} is exact so $\im(d_{n+1})=\ker(d_{n})$ for $n\geq 1$. Finally $\im(d^\mu_1)=F_0(\im(d_1))$, but $F_0=\id$, so $\im(d^\mu_1)=\im(d_1)=\ker(\epsilon')$, again using the fact \eqref{augmented_complex} is exact. Therefore $K_\bullet(A_\mu)$ is a resolution of $k$ as an $A_\mu$-module, and $A_\mu$ is a Koszul algebra.

\section{Examples}\label{examples_sec}
Next we demonstrate Theorems \ref{quadratic_h_mod_alg_thm} and \ref{koszul_main_thm} with several examples. In particular, in Section \ref{plane_sec}, we use Theorem \ref{quadratic_h_mod_alg_thm} to determine the result of twisting the quantum plane by the quasitriangular structure of the quantum enveloping algebra $U_q(\mathfrak{sl}_2)$, whilst in Section \ref{sym_ext_sec} we apply Theorem \ref{koszul_main_thm} to provide a new proof of the Koszulity of the quantum symmetric and exterior algebras $S_{-1}(V)$ and $\bigwedge_{-1}(V)$.

\subsection{The quantum plane and \texorpdfstring{$U_q(\mathfrak{
sl}_2)$}{sl}.}\label{plane_sec}
For a non-zero $q\in k$, the $q$-quantum plane is given by $A=k\langle x,y\rangle/(xy-qyx)$. In the following  we explain why we can twist $A$ by the quasitriangular structure $\mathcal{R}$ of the quantum enveloping algebra $U_q(\mathfrak{sl}_2)$, and we find via Theorem \ref{quadratic_h_mod_alg_thm} that $A_\mathcal{R}$ is the quadratic algebra equal to the $q^{-1}$-quantum plane $k\langle x,y\rangle/(xy-q^{-1}yx)$.\bb

\nt First recall that a quasitriangular structure on a Hopf algebra $H$ is an invertible $\mathcal{R}\in H\otimes H$ satisfying: 
\begin{equation}\label{qt_structure}
\tau\circ \triangle(h)=\mathcal{R}\cdot \triangle(h)\cdot \mathcal{R}^{-1},\ \forall h\in H,\hspace{1em} (\triangle\otimes \id)(\mathcal{R})=\cR_{13}\cdot \cR_{23},\hspace{1em} (\id\otimes \triangle)(\cR)=\cR_{13}\cdot \cR_{12}
\end{equation}
where $\tau$ is the transposition map $v\otimes w\mapsto w\otimes v,\ \cR_{12}=\cR\otimes 1,\ \cR_{23}=1\otimes R$ and $\cR_{13}$ has $1$ inserted in the middle leg. Quasitriangular structures are known to satisfy the quantum Yang-Baxter equation \cite[Lemma 2.1.4]{alma998944944401631}: $\cR_{12}\cdot \cR_{13}\cdot \cR_{23}=\cR_{23}\cdot \cR_{13}\cdot \cR_{12}$, and from this, one can check that a quasitriangular structure is also a counital $2$-cocycle, i.e. satisfies \eqref{2_cocycle_def} and \eqref{counital_def}. For this reason it makes sense to use quasitriangular structures to perform Drinfeld twists.\bb

\nt Let us introduce the quantum enveloping algebra $U_q(\mathfrak{sl}_2)$. We now suppose $q$ also satisfies $q^2\neq 1$. Then $U_q(\mathfrak{sl}_2)$ is defined to be the $k$-algebra generated by $E,F, K,K^{-1}$, and satisfying the relations
$$KE=q^2 EK,\hspace{1em} KF=q^{-2}KF,\hspace{1em} EF-FE=\frac{K-K^{-1}}{q-q^{-1}},\hspace{1em} KK^{-1}=K^{-1}K=1.$$
\nt $U_q(\mathfrak{sl}_2)$ may additionally be given the structure of a Hopf algebra via the following:
\begin{align*}
\triangle(E) & = 1\otimes E+E\otimes K,& \epsilon(E) & =0,& S(E) & =-EK^{-1}\\ 
\triangle(F) & =K^{-1}\otimes F + F\otimes 1,& \epsilon(F) & =0,& S(F) & =-KF\\
\triangle(K^{\pm 1}) & =K^{\pm 1}\otimes K^{\pm 1},& \epsilon(K^{\pm 1}) & =1,& S(K^{\pm 1}) & =K^{\mp 1}
\end{align*}
Now the quantum plane $A=k\langle x,y\rangle/(xy-qyx)$ is, by construction, a quadratic algebra. Additionally, we may define a representation of $U_q(\mathfrak{sl}_2)$ on the degree $1$ homogeneous subspace $V=\spn_k\{x,y\}$ of $A$ as follows:
\begin{align*}
E\rhd x & =0,& F\rhd x& =y,& K^{\pm 1}\rhd x=q^{\pm 1}x,\\
E\rhd y & =x,& F\rhd y& =0,& K^{\pm 1}\rhd y=q^{\mp 1}y
\end{align*}
It is well-known (see, for instance, \cite[Exercise 9.1.13]{alma998944944401631}) that this action extends to make $A$ a $U_q(\mathfrak{sl}_2)$-module algebra. Note that, by construction, this action is degree-preserving.\bb

\nt In the terminology of Vlaar \cite[Theorem 6.7]{vlaar2020}, $U_q(\mathfrak{sl}_2)$ has a quasitriangular structure $\cR$ ``up to completion" - meaning that $\cR$ does not lie in $U_q(\mathfrak{sl}_2)\otimes U_q(\mathfrak{sl}_2)$, but rather in a completion of this space. Despite this technicality, $\mathcal{R}$ still satisfies axioms \eqref{2_cocycle_def} and \eqref{counital_def} of a counital $2$-cocycle. However, since $\cR$ is not in $U_q(\mathfrak{sl}_2)\otimes U_q(\mathfrak{sl}_2)$, we must check that there is still a well-defined action of $\cR$ on $A\otimes A$ in order for us to define the twist of $A$ by $\cR$.\bb 

\nt Etingof \cite[Remark 3.41]{etingof2021brief} states that given two representations of $U_q(\mathfrak{sl}_2)$, say $\rho:U_q(\mathfrak{sl}_2)\rt \End(V)$ and $\rho':U_q(\mathfrak{sl}_2)\rt \End(W)$, which are locally nilpotent (i.e. $\forall v\in V$ or $W$, $\exists n\in \N$ such that $E^n\rhd v=0$), we have that $(\rho\otimes \rho')(\cR)$ is a well-defined operator on $V\otimes W$. Therefore, if the action of $U_q(\mathfrak{sl}_2)$ on $A$ is locally nilpotent then it will follow that $\cR$ has a well-defined action on $A\otimes A$. The action on $A$ is indeed locally nilpotent, and this can be seen from how $E$ acts on a general basis element of $A$:
\begin{equation}\label{uqsl2_action}
E\rhd x^a y^b=[b]_qx^{a+1}y^{b-1},\hspace{2em}\text{where }[b]_q:=\frac{q^b-q^{-b}}{q-q^{-1}}
\end{equation}
It is a simple exercise to check \eqref{uqsl2_action}, first by showing $E\rhd x^a y^b=x^a (E\rhd y^b)$, and then using induction (in the degree $b$), to show $E\rhd y^b=[b]_q xy^{b-1}$. From \eqref{uqsl2_action}, we see $E^{b+1}\rhd x^ay^b=0$, and so $U_q(\mathfrak{sl}_2)$ indeed acts locally nilpotently on $A$. Hence $\cR$ has a well-defined action on $A\otimes A$ and 
we may construct the Drinfeld twist $A_\cR$.\bb 

\nt Finally we can apply our first main result, Theorem \ref{quadratic_h_mod_alg_thm}. The conditions of the theorem are met since the action of $U_q(\mathfrak{sl}_2)$ on $A$ is degree-preserving. We deduce that $A_\cR$ is a quadratic algebra, and is given by $k\langle x,y\rangle / (\cR \rhd (x\otimes y-q y\otimes x))$. Vlaar \cite[Equation 6.37]{vlaar2020} tells us how $\cR$ acts on $V\otimes V$ explicity, and for the basis vectors $x\otimes y$ and $y\otimes x$ it is,
$$\cR\rhd x\otimes y=q^{-1/2}(x\otimes y+(q-q^{-1})y\otimes x),\hspace{2em}\cR\rhd y\otimes x=q^{-1/2}y\otimes x$$
Therefore 
$$\cR \rhd x\otimes y-q y\otimes x=q^{-1/2}(x\otimes y-q^{-1}y\otimes x)$$
and so $\cR \rhd (x\otimes y-q y\otimes x)$ is equal to the ideal $(x\otimes y-q^{-1}y\otimes x)$. Hence $A_\cR$ is the $q^{-1}$-quantum plane.

\subsection{Symmetric and Exterior algebras.}\label{sym_ext_sec}
Here we will consider $k=\Cc$, and in each example we will take the Hopf algebra $H$ to be the group algebra $\Cc T$, where $T$ is the finite abelian group given by $(C_2)^n=\langle t_1,\dots,t_n\ |\ t^2_i=1, t_it_j=t_jt_i\rangle$ for some $n\geq 2$. The counit of $H$ is given by $\epsilon(t)=1,\ \forall t\in T$, whilst the coproduct is $\triangle(t)=t\otimes t\ \forall t\in T$.\bb

\nt Let, 
$$\mu=\prod_{1\leq j<i\leq n}\mu_{ij},\hspace{2em}\mu_{ij}=\frac{1}{2}(1\otimes 1+t_i\otimes 1+1\otimes t_j-t_i\otimes t_j)$$
By \cite[Lemma 4.5]{twistsrcas}, $\mu$ is a counital $2$-cocycle of $H=\Cc T$.  Note also that $\mu=\mu^{-1}$.\bb

\nt Consider an $n$-dimensional $\Cc $-vector space $V$ with a fixed basis $x_1,\dots,x_n$. For an $n\times n$-matrix $q=(q_{ij})$ satisfying $q_{ii}=1$ and $q_{ij}q_{ji}=1$, the corresponding \define{quantum symmetric algebra} is defined as 
$$S_q(V):=T(V)/(x_i\otimes x_j-q_{ij}x_j\otimes x_i\ |\ 1\leq i,j \leq n).$$
Likewise the \define{quantum exterior algebra} is given as $\bigwedge_q(V):=T(V)/ (x_i \otimes x_j+q_{ij} x_j\otimes x_i\ |\ 1\leq i,j \leq n)$. In the following examples we will be interested in the case when $q=(-1)$, where $$(-1)_{ij}=\begin{cases}
1, & \text{if i=j}\\
-1, & \text{otherwise.}
\end{cases}$$
The quantum symmetric and exterior algebras are known to be Koszul, but we show that this can also be deduced as an application of Theorem \ref{koszul_main_thm}.

\begin{example}[Twisting the symmetric algebra.]
Recall that the symmetric algebra $S(V)\cong \Cc[x_1,\dots,x_n]$ is a Koszul algebra (see Witherspoon \cite[Example 3.4.11]{alma992981689925201631}), and by definition is equal to $T(V)/(R)$ for $R=\spn_\Cc\{x_i\otimes x_j-x_j\otimes x_i\ |\ 1\leq i,j\leq n\}$. Additionally $H=\Cc T$ acts on $V$ via $t_i\rhd x_j=(-1)^{\delta_{ij}}x_j$, and this extends by algebra homomorphisms to make $S(V)$ an $H$-module algebra. In particular the action on a monomial is given by $$t_i\rhd x_1^{a_1}\dots x_i^{a_i}\dots x_n^{a_n}=(-1)^{a_i}x_1^{a_1}\dots x_i^{a_i}\dots x_n^{a_n}$$
Note that this action of $H$ on $S(V)$ is degree-preserving, and therefore we may apply Theorem \ref{koszul_main_thm} to deduce that the Drinfeld twist $S(V)_\mu$ is a Koszul algebra isomorphic to $T(V)/(R_\mu)$, where 
$$R_\mu=\spn_\Cc \{\mu\rhd (x_i\otimes x_j-x_j\otimes x_i)\ |\ 1\leq i,j\leq n\}$$
By \cite[Corollary 5.8]{twistsrcas}, for $i\neq j$, $\mu\rhd (x_i\otimes x_j-x_j\otimes x_i) =x_i\otimes x_j+x_j\otimes x_i$, so $S(V)_\mu=S_{-1}(V)$. Therefore we find that Theorem \ref{koszul_main_thm} provides a new proof that the quantum symmetric algebra $S_{-1}(V)$ is Koszul.
\end{example}

\begin{example}[Twisting the exterior algebra.]
The exterior algebra is also Koszul, and is equal to $T(V)/(R)$ for $R=\spn_\Cc\{x_i\otimes x_j+x_j\otimes x_i\ |\ 1\leq i,j\leq n\}$. The action of $H=\Cc T$ on $V$ given by $t_i\rhd x_j=(-1)^{\delta_{ij}}x_j$ also extends by algebra automorphisms to the exterior algebra $\bigwedge (V)$, making $\bigwedge(V)$ an $H$-module algebra. By construction this action is degree-preserving, and so, by Theorem \ref{koszul_main_thm}, we find $\bigwedge(V)_\mu$ is a Koszul algebra. Now  $\bigwedge(V)_\mu=T(V)/(R_\mu)$ for
\begin{align*}
R_\mu & =\spn_\Cc \{\mu\rhd (x_i\otimes x_j+x_j\otimes x_i)\ |\ 1\leq i,j\leq n\}\\
& = \spn_\Cc \{x_i\otimes x_i,\ x_i\otimes x_j-x_j\otimes x_i\ |\ i\neq j\}
\end{align*}
where we apply \cite[Corollary 5.8]{twistsrcas} in the second equality. Therefore $\bigwedge(V)_\mu=\bigwedge_{-1}(V)$, and we have found another proof that the quantum exterior algebra $\bigwedge_{-1}(V)$ is Koszul.
\end{example}
\bibliographystyle{plain}
\bibliography{mybib}

\begin{thebibliography}{1}

\bibitem{twistsrcas}
Y.~Bazlov, A.~Berenstein, E.~Jones-Healey, and A.~Mcgaw.
\newblock {Twists of rational Cherednik algebras}.
\newblock {\em The Quarterly Journal of Mathematics}, 74(2):511--528, 2022.

\bibitem{doi:10.1080/00927872.2016.1178271}
A.~Davies.
\newblock Cocycle twists of {A}lgebras.
\newblock {\em Communications in Algebra}, 45(3):1347--1363, 2017.

\bibitem{drinfeld_quantum_groups}
V.~G. Drinfel'd.
\newblock Quantum {G}roups.
\newblock {\em Proceedings of the International Congress of Mathematicians},
  Vol. 1:798--820, 1986.

\bibitem{etingof2016tensor}
P.~Etingof, S.~Gelaki, D.~Nikshych, and V.~Ostrik.
\newblock {\em Tensor Categories}, volume 205 of {\em Mathematical Surveys and
  Monographs}.
\newblock American Mathematical Society, Providence, RI, 2016.

\bibitem{etingof2021brief}
P.~Etingof and M.~Semenyakin.
\newblock A {B}rief {I}ntroduction to {Q}uantum {G}roups.
\newblock {\em arXiv:2106.05252}, 2021.

\bibitem{giaquinto1998bialgebra}
A.~Giaquinto and J.~J. Zhang.
\newblock Bialgebra actions, twists, and universal deformation formulas.
\newblock {\em Journal of Pure and Applied Algebra}, 128(2):133--151, 1998.

\bibitem{alma998944944401631}
S.~Majid.
\newblock {\em Foundations of Quantum Group Theory}.
\newblock Cambridge University Press, Cambridge, UK, 1995.

\bibitem{vlaar2020}
B.~Vlaar.
\newblock {LMS} {A}utumn {A}lgebra {S}chool 2020, {L}ecture notes: Introduction
  to {Q}uantum {G}roups.
\newblock
  \url{https://www.icms.org.uk/sites/default/files/documents/events/Bart%20Vlaar.pdf},
  2020.
\newblock [Online; accessed 31-May-2023].

\bibitem{alma992981689925201631}
S.~J. Witherspoon.
\newblock {\em Hochschild {C}ohomology for {A}lgebras}, volume 204 of {\em
  Graduate {S}tudies in {M}athematics}.
\newblock American Mathematical Society, Providence, RI, 2019.

\end{thebibliography}

\end{document}